\documentstyle[12pt,aasms4]{article}

\renewcommand{\baselinestretch}{2}

\newcommand{\SS}{\renewcommand{\baselinestretch}{1}   \tiny \noindent}

\begin{document}

\SS

\section*{ }

\begin{center}
\Huge
{\bf  Bayesian Estimation of\\
  Time Series Lags and Structure}
\vskip 0.5in
\Large
{\bf Jeffrey D. Scargle}\\
{\bf Space Science Division}\\
{\bf NASA Ames Research Center}\\
{\bf MS 245-3, NASA Ames Research Center}\\
{\bf Moffett Field, CA, 94035-1000}\\
{\bf \verb+jeffrey@cosmic.arc.nasa.gov+}\\

\end{center}

% \vskip 1.75in
% \noindent
% Contribution to {\bf Workshop on Bayesian Inference and Maximum Entropy Methods in
% Science and Engineering (MAXENT 2001)}, held at Johns Hopkins University, 
% Baltimore, MD USA on August 4-9, 2001.

\large

%================================================

% \tableofcontents

\section*{Abstract}
This paper derives practical algorithms, based on
Bayesian inference methods, for several data 
analysis problems common in time series analysis
of astronomical and other data.  One problem
is the determination of the lag between two time
series, for which the {\it cross-correlation function} 
is a sufficient statistic.
The second problem is the estimation of structure in a time
series of measurements which are a weighted 
integral over a finite range of the independent variable.  

% \newpage

\section{Workhorse Algorithms for Time Series Analysis}

Bayesian methods are becoming more popular
for the challenging data analysis problems facing the
modern astrophysicist, but the pace is agonizingly slow.
I believe the main difficulty is the perception that Bayesian
methods must be implemented in complex, special-purpose routines
made for a single application and requiring copious computational resources.
Progress will be accelerated by the availability of
turn-key algorithms for elementary data analysis problems.

Larry Bretthorst has pioneered in developing Bayesian
methods for the detection of periodic signals in noisy data
(Bretthorst 1988, 2001).
He computes the posterior distribution of the
frequency parameter in a model consisting of 
a single sinusoidal component, having marginalized the
amplitude and phase parameters.
It is encouraging that this work is making its way into 
a number of active areas in astronomy,
including variable star research and more recently 
discovery of extra-solar planets.
The present work applies the methods clearly
outlined in (Bretthorst 1988) to another common
astronomical problem -- the detection of lags between
two or more signals.

\section{Lags in Time Series Data}

In engineering and science, including both
experimental and observational sciences, such as astronomy, 
one often wishes to find the delay between two time series.  
This somewhat complex mixture of questions includes:
Are the two time series related?  
If they are, is one a delayed version of the other?  
If so, what is the best estimate of the value of the lag?

\subsection{The Model}

A straightforward approach is to define 
a generic model expressing one signal as a delayed and scaled 
version of the other, and then derive the posterior 
probability distribution of the parameters
representing the lag and the scale factor.
This procedure can be carried out making few
assumptions about the signal, and none about signal shape. 
From this posterior one can easily compute 
means and confidence intervals for lags and scale factors.

Let $X$ and $Y$ denote two observables.
Assume that the underlying process being sampled here
is an unknown signal, $S$, superimposed on a constant background, $B$.
If a negative signal is impossible for physical reasons
an appropriate prior can impose the condition $S \ge 0$.

The backgrounds can be treated as unknown nuisance parameters, 
assigned a prior probability distribution, and marginalized.  
If the backgrounds are accurately fixed by other data, 
so the prior distribution is very narrow, one can sometimes
get away with treating the backgrounds as known constants.

The model of the observables, expressing delay and scaling
between the two signals, is then:
\begin{eqnarray} \label{y_model_1}
X_{m} = &  S_{m}         & + \ B_{X} \\  
Y_{m} = & a S_{m - \tau} & + \ B_{Y} \ , \label{y_model_2}
\end{eqnarray}
\noindent
where $m$ represents the independent variable, often time,
the lag is $\tau$, and we allow 
the $Y$-signal to be an overall factor
$a$ times the $X$-signal.
Of course, $a$ may be less than, equal to,
or greater than $1$.

We now discuss two data modes common in astronomy,
namely time-tagged events and 
evenly sampled time series with normal errors.
The different nature of the observational errors 
in these two cases means that they are 
represented differently in the model,
as will be seen in the next two sections.

\subsection{Time-Tagged Event (TTE) Data}

We begin by treating event data,
sometimes called {\it time-tagged event} (TTE) data in
the astronomical literature.
Such data are simply the set of times 
at which events occurred -- usually 
within a fixed interval, 
starting at time $0$ and ending
at time $T$.
Here we assume that the only observational
noise is due to the randomness of the events.\footnote{E.g.,
photon detection, the most common astronomical application,
is inherently discrete due to the quantum nature of light.
The physical parameter is the expected rate of photon
detection, determined by the incident photon intensity
and the instrument's detection efficiency.}

The times are not, of course, recorded with infinite precision.
They are quantized in small units, 
here called {\it time-ticks}, defined by the
computer clock that drives the data acquisition system.
Setting the time-tick to unity,
the event times are a set of integers satisfying
\begin{equation}
0 \le m_{1} < m_{2} < m_{3} < \dots < m_{N-1} < m_{N} \le T
\label{data}
\end{equation}
\noindent
where $N$ is the total number of events.
Often the detection process
mandates the condition $m_{i} \ne m_{i+1}$ indicated
in Eq. (\ref{data}).\footnote{
In a few cases -- {\it e.g.} multiple detectors 
in a single spacecraft -- this is not true.}
Indeed, it is almost always the case
that each event is followed by a short interval 
during which the instrument cannot detect any
subsequent event.  We ignore this {\it detector dead time}.

It is useful to represent TTE data
as a series of $N$ delta functions:
\begin{equation}
x_{m} = \left\{ \begin{array}{llll}
&1 \  & \mbox{if event at time $m$ } & \\
&  \  &  & m = 1, 2, \dots, M \\
&0 \  & \mbox{if no event at time $m$} & 
\end{array}
\right.
\label{data_1}
\end{equation}
\noindent
where $M$ is the total length of the observation
interval in time-ticks,
and $x$ is the observed value of $X$
(similarly for $y$ and $Y$).

As mentioned above, 
we assume that the only
observational noise is that due to the randomness of the events.
Equations (\ref{y_model_1}) and (\ref{y_model_2}) give
the probability of detecting an event 
during tick $m$.
Typically the instrument is designed 
so that these probabilities $X_{m}, Y_{m}$ are $<<1$.
In this {\it truncated Poisson process}
the probability of no $X$-event at time $m$ is $e^{-X_{m}}$.
Hence the likelihood is simply
\begin{equation}
L( x_{m} | S_{m}, B_{X} ) = \left\{ \begin{array}{lc}
e^{- (S_{m} + B_{X}) } & x_{m} = 0 \\
1 - e^{ - ( S_{m} + B_{X} )}  & x_{m} = 1
\end{array}
\right.
\label{likelihood_x}
\end{equation}
Similarly for $Y$
\begin{equation}
L( y_{m} | S_{m-\tau},  B_{Y}, a ) = \left\{ \begin{array}{lc}
e^{- ( a S_{m-\tau} + B_{Y}) } & y_{m} = 0 \\
1 - e^{ - ( a S_{m-\tau} + B_{Y} )}  & y_{m} = 1
\end{array}
\right.
\label{likelihood_y}
\end{equation}
\noindent
It is fundamental to this analysis that the
$X_{m}$ and $Y_{m}$ are all independent 
with respect to the measurement noise process.
Hence the total likelihood is 
the product of the individual ones:
\begin{equation}
L_{total} = \prod_{m=1}^{M} 
L( x_{m} | S_{m},         B_{X} ) 
L( y_{m} | S_{m - \tau},  B_{Y}, a )
\label{total_like_1}
\end{equation}
Note that with these likelihoods 
there is an issue connected with
{\it wraparound} that is essentially the same
as with cross correlation functions of any kind.
The expressions derived here 
assume that the data may be allowed
to wraparound.

It is more convenient to use the form
\begin{equation}
Y_{m + \tau} = a S_{m} + B_{Y} \ ,
\end{equation}
equivalent to Eq. (\ref{y_model_2}),
to transform the expression for $L$ in Eq. (\ref{total_like_1}) to
\begin{equation}
L_{total} = \prod_{m=1}^{M} 
L( x_{m}        | S_{m}, B_{X}  ) 
L( y_{m + \tau} | S_{m}, B_{Y} , a )
\label{total_like_2}
\end{equation}
so that we can write the total likelihood
as the product of factors, each of which
depends on the same signal variable, $S_{m}$:
\begin{equation}
L_{total} = \prod_{m=1}^{M} L_{m}(S_{m})
\label{total_like_3}
\end{equation}
\noindent
where
\begin{equation}
L_{m}(S_{m}) =  
L( x_{m}      | S_{m}, B_{X} ) 
L( y_{m+\tau} | S_{m}, B_{Y}, a )
\label{like_m}
\end{equation}

We can individually marginalize the signal
parameters $S_{m}$, which for
the purposes of determining the
lag and scale factor 
are nuisance parameters.
Dropping the subscript on the now dummy varible
$S_{m}$,
and adopting a prior $P( S )$, the marginalized posterior is 
\begin{equation}
G(x_{m},y_{m+\tau}) =  \int _{-\infty}^{+\infty}
P( S )
L( x_{m}      | S, B_{X} )
L( y_{m+\tau} | S, B_{Y}, a ) dS
\label{post_m}
\end{equation}
That is to say, we have
\begin{equation}
G_{total}(\tau, a)
= \prod_{m=1}^{M} G( x_{m}, y_{m+\tau} )
\end{equation}
\noindent
or
\begin{equation}
G_{total}(\tau, a)
= G_{0,0}^{N_{0,0}} G_{0,1}^{N_{0,1}} G_{1,0}^{N_{1,0}} G_{1,1}^{N_{1,1}}
\label{posterior_factors}
\end{equation}
where 
\begin{eqnarray}
N_{0,0}(\tau) = & \mbox{number of \ } m \ \mbox{\ for which} & x_{m} = 0 \ \mbox{and} \ y_{m+ \tau} = 0 \\
N_{0,1}(\tau) = &                            & x_{m} = 0 \ \mbox{and} \ y_{m+ \tau} = 1 \\
N_{1,0}(\tau) = &                            & x_{m} = 1 \ \mbox{and} \ y_{m+ \tau} = 0 \\
N_{1,1}(\tau) = &                            & x_{m} = 1 \ \mbox{and} \ y_{m+ \tau} = 1 
\end{eqnarray}
and
\begin{equation}
\begin{array}{cc}
G_{\alpha,\beta} = \int
P( S )
L( x_{m} = \alpha     | S, B_{X} ) 
L( y_{m+\tau} = \beta | S, B_{X}, a ) dS \\ 
 \ \ \\
\alpha = 0,1; \ \beta = 0,1
\end{array}
\end{equation}
\noindent
Note that 
the $N$'s
depend on the data and on $\tau$, but not on $a$; 
the $G$'s depend on 
the model's form,
priors on model parameters,
and on $a$, but not on $\tau$.

The $N$'s as functions of lag $\tau$ 
are conveniently found from
the cross-correlation
function of $X$ and $Y$,
defined as
\begin{equation}
\gamma_{X,Y}(\tau) = \sum_{m=1}^{M}
 X_{m+\tau} y_{m}  
\label{cross_deff}
\end{equation}
This function is readily and rapidly
computed, using the fast Fourier transform,
by representing 
$X$ and $Y$
as arrays of zeros punctuated by unit amplitude
$\delta$-functions at the $m$ at which events occur.

It can be shown that 
\begin{eqnarray}
N_{1,1} = & \gamma_{X,Y}(\tau) \label{cc1} \\
N_{1,0} = & N_{X} - \gamma_{X,Y}(\tau) \label{cc2} \\ 
N_{0,1} = & N_{Y} - \gamma_{X,Y}(\tau) \label{cc3} \\ 
N_{0,0} = & M - N_{X} - N_{Y} + \gamma_{X,Y}(\tau) \label{cc4}
\end{eqnarray}
where $M$ is the number of values of $m$ spanning the observation interval
[{\it cf.} Eq. (\ref{data_1})],
and $N_{X}$ and $N_{Y}$ 
are just the number of events.
Since the four combinations
exhaust all possibilities,
these quantities should, and obviously do, satisfy
\begin{equation}
N_{0,0} + N_{0,1} + N_{1,0} + N_{1,1} = M
\end{equation}

Adopting the uniform prior 
\begin{equation}
P( S ) = {1 \over S_{1} - S_{0} }
\left\{ \begin{array}{cl}
   1 & \mbox{for $ S_{0} \le S \le S_{1} $ } \\
   0 & \mbox{else}
\end{array}
\right .
\end{equation}
the $G$'s are easily found to be:
\begin{eqnarray}
G_{1,1} = & 1 - E_{X}\phi( S_{0}, S_{1} ) 
- E_{Y}\phi( a S_{0}, a S_{1} )
+ E_{X}E_{Y} \phi( \rho S_{0}, \rho S_{1} ) \\    \label{n1}
G_{1,0} = &  E_{Y}\phi( aS_{0}, aS_{1} ) -
 E_{X}E_{Y}\phi( \rho S_{0}, \rho S_{1} ) \\    \label{n2}
G_{0,1} = & E_{X}\phi( S_{0}, S_{1} ) -
 E_{X}E_{Y}\phi( \rho S_{0}, \rho S_{1} ) \\    \label{n3}
G_{0,0} = & E_{X} E_{Y} \phi( \rho S_{0}, \rho S_{1} )  \label{n4}
\end{eqnarray}
\noindent
where $\rho = 1 + a$,
\begin{equation}
E_{X} = e^{-B_{X}}
\end{equation}
\noindent
\begin{equation}
E_{Y} = e^{-B_{Y}}
\end{equation}
\noindent
and
\begin{equation}
\phi( x, y ) \equiv {e^{-x} - e^{-y} \over x - y}
\end{equation}

It is instructive to take log 
of the likelihoods in Eq.(\ref{posterior_factors}), as follows:
\begin{equation}
log G_{total}(\tau,a)
\equiv
N_{0,0} log G_{0,0} + N_{0,1} log G_{0,1} + N_{1,0} log G_{1,0} + N_{1,1} log G_{1,1}
\end{equation}
and Equations (\ref{cc1}-\ref{cc4})
permit a representation in the form
\begin{equation}
log G_{total}(\tau,a) = c_{0}(a) + c_{1}(a) \ \gamma_{X,Y}(\tau)
\end{equation}
where 
\begin{equation}
c_{0}(a) = [M - N_{X} - N_{Y}] log G_{0,0} 
+ N_{Y} log G_{0,1} +  N_{X} log G_{1,0}
\end{equation}
and
\begin{equation}
c_{1}(a) = log G_{0,0} - log G_{0,1} - log G_{1,0} + log G_{1,1}
\end{equation}

Accordingly, we have for the 
posterior probability density
\begin{equation}
G_{total}(\tau, a) = e^{c_{0}(a)} e^{ c_{1}(a) \ \gamma_{X,Y}(\tau)  }
\propto e^{ c_{1}(a) \ \gamma_{X,Y}(\tau)  } 
\label{posterior}
\end{equation}
Note the similarity of Eq. (\ref{posterior})
to an analogous result in harmonic analysis -- detecting a sinusoidal
signal in the presence of noise --
giving the posterior probability density for the frequency $\omega$ 
[\cite{bretthorst_1}, Eq. (2.7)] 
\begin{equation}
P( \omega | D, \sigma, I) \propto e^{ {C(\omega) \over \sigma^{2} } }
\label{bretthorst}
\end{equation}
where
$C(\omega)$ is the periodogram, 
$D$ represents the data, $\sigma$ is the variance
of the noise, here assumed known, and
$I$ is the prior information.
The cross correlation function, $\gamma$
is a sufficient statistic for lags,
just as the periodogram is for frequencies (\cite{bretthorst_1}).

Note that the maximum likelihood 
value of the lag $\tau$ is just the
value that maximizes the cross-correlation function.
The mean value of $\tau$, weighted 
by the posterior in Eq. (\ref{posterior}), may be a better
lag estimator.  This posterior is also useful for
computing confidence intervals.

\subsection{Evenly Spaced Data}

For simplicity, in this section we
ignore the background component, and
consider noisy measurements of two signals
\begin{eqnarray}
X_{m} & = S^{X}_{m} &+ \ R^{X}_{m}\\
Y_{m} & = S^{Y}_{m} &+ \ R^{Y}_{m}
\label{even_data}
\end{eqnarray}
where $m = 1, 2, \dots, M$ are 
a set of evenly spaced times,
and the true signals $S^{X,Y}$ 
are corrupted by additive noise, $R^{X,Y}$.
Consider the model
\begin{equation}
S^{Y}_{m} = a S^{X}_{m-\tau}
\label{model_sampled}
\end{equation}
stating that one signal is a delayed and
scaled version of the other.

Assume that the noise has
a normal distribution
[see (\cite{bretthorst_1}) for
relevant discussion], so that the 
likelihood for $X$ at time $m$ is
\begin{equation}
P( x_{m} | S^{X}_{m}, \sigma_{m}^{X} ) = 
{1 \over \sigma_{m}^{X} \sqrt{ 2 \pi  } }
e^{ - { R_{m}^{2} \over 2 (\sigma_{m}^{X})^{2} } } =
{1 \over \sigma_{m}^{X} \sqrt{ 2 \pi  } }
e^{ - { ( x_{m} - S^{X}_{m} )^{2} \over 2 (\sigma_{m}^{X})^{2} } }
\end{equation}
Similarly for $Y$ at time $m$
\begin{equation}
P( y_{m} | S^{Y}_{m}, \sigma_{m}^{Y} ) = 
{1 \over \sigma_{m}^{Y} \sqrt{ 2 \pi  } }
e^{ - { ( y_{m} - S^{Y}_{m} )^{2} \over 2 (\sigma_{m}^{Y})^{2} } }
\end{equation}
\noindent
which, from Eq. (\ref{model_sampled}), becomes
\begin{equation}
P( y_{m} | S^{X}_{m-\tau}, a, \sigma_{m}^{Y} ) =
{1 \over \sigma_{m}^{Y} \sqrt{ 2 \pi  } }
e^{ - { ( y_{m} - aS^{X}_{m-\tau} )^{2} \over 2 (\sigma_{m}^{Y})^{2} } }
\end{equation}
\noindent
With the usual independence assumption, the total likelihood is
\begin{equation}
P_{total} = ({1 \over 2 \pi })^{M}
\prod_{m=1}^{M}
{1 \over \sigma_{m}^{X} \sigma_{m}^{Y} }
e^{ - { ( x_{m} - S^{X}_{m} )^{2} \over 2 (\sigma_{m}^{X})^{2} } 
    - { ( y_{m} - aS^{X}_{m-\tau} )^{2} \over 2 (\sigma_{m}^{Y})^{2} } }
\end{equation}
That is 
\begin{equation}
P_{total} = P_{0}
\prod_{m=1}^{M}
e^{ - { ( x_{m} - S^{X}_{m} )^{2} \over 2 (\sigma_{m}^{X})^{2} } 
    - { ( y_{m} - aS^{X}_{m-\tau} )^{2} \over 2 (\sigma_{m}^{Y})^{2} } }
\end{equation}
where
\begin{equation}
P_{0} = ({1 \over 2 \pi })^{M}
\prod_{m=1}^{M}
{1 \over \sigma_{m}^{X} \sigma_{m}^{Y} }
\end{equation}

Shifting $m \rightarrow m+\tau$ in part of this expression
gives an equivalent form in which factors involving
the same $S_{m}$ are kept together:
\begin{equation}
P_{total} = P_{0}
\prod_{m=1}^{M}
e^{ - { ( x_{m} - S^{X}_{m} )^{2} \over 2 (\sigma_{m}^{X})^{2} } 
    - { ( y_{m+\tau} - aS^{X}_{m} )^{2} \over 2 (\sigma_{m+\tau}^{Y})^{2} } }
\equiv P_{0}\prod_{m=1}^{M}  e^{L_{m}}
\end{equation}
Expanding the argument of the exponential:
\begin{equation}
L_{m} \equiv
- { x_{m}^{2} - 2 x_{m} S^{X}_{m} + (S^{X}_{m})^{2} \over 2 (\sigma_{m}^{X})^{2} } 
- { y_{m+\tau}^{2} - 2aS^{X}_{m}y_{m+\tau} + (aS^{X}_{m}) ^{2} \over 2 (\sigma_{m+\tau}^{Y})^{2} } 
\end{equation}
\noindent
or
\begin{equation} =
- {1 \over 2}[ { x_{m}^{2} \over (\sigma_{m}^{X})^{2} }
 + { y_{m+\tau}^{2} \over  (\sigma_{m+\tau}^{Y})^{2} } ]
\ \ 
+ [{ x_{m}  \over  (\sigma_{m}^{X})^{2} } +
{ a y_{m+\tau} \over  (\sigma_{m+\tau}^{Y})^{2} } ] S^{X}_{m}
\ \ 
- {1 \over 2} [ { 1 \over (\sigma_{m}^{X})^{2} }
 \ + { a^{2} \over (\sigma_{m+\tau}^{Y})^{2} } ] (S^{X}_{m})^{2}
\end{equation}
Rewrite this as
\begin{equation}
L_{m} = -C_{m} + 2 B_{m} S^{X}_{m}  - A_{m}(S^{X}_{m})^{2}
\end{equation}
\noindent
where
\begin{equation}
C_{m} \equiv
{1 \over 2}[ { x_{m}^{2} \over (\sigma_{m}^{X})^{2} }
 + { y_{m+\tau}^{2} \over  (\sigma_{m+\tau}^{Y})^{2} } ]
\end{equation}
\begin{equation}
B_{m} \equiv
{1 \over 2} [{ x_{m}  \over  (\sigma_{m}^{X})^{2} } +
{ a y_{m+\tau} \over  (\sigma_{m+\tau}^{Y})^{2} } ]
\end{equation}
\begin{equation}
A_{m} \equiv
{1 \over 2} [ { 1 \over (\sigma_{m}^{X})^{2} }
 \ + { a^{2} \over (\sigma_{m+\tau}^{Y})^{2} } ] 
\end{equation}
Still following (\cite{bretthorst_1}) we complete the square:
\begin{equation}
L_{m} = - A_{m} ( S^{X}_{m} - {B_{m} \over A_{m} })^{2} - C_{m}
+ {B_{m}^{2} \over A_{m}}
\end{equation}
\noindent
so that the marginalizations of the $S_{m}$ become
\begin{eqnarray}
P( \tau, a | D ) & = & P_{0}
\prod_{m=1}^{M} \int_{-\infty}^{\infty} e^{ -A_{m} ( S^{X}_{m} - {B_{m} \over A_{m} })^{2}
- C_{m}
+ {B_{m}^{2} \over A_{m}}} \ \  dS_{m} \\
\ & = & P_{0}
\prod_{m=1}^{M} e^{ {B_{m}^{2} \over A_{m}} - C_{m} }
\int_{-\infty}^{\infty} e^{ -A_{m} ( S^{X}_{m} - {B_{m} \over A_{m} })^{2}} dS_{m} \\
\ & = & P_{0}
\prod_{m=1}^{M} e^{ {B_{m}^{2} \over A_{m}} - C_{m} }
\sqrt{ {\pi \over A_{m} } } \\
 & = & P_{0}
e^{ \sum_{m=1}^{M} ( {B_{m}^{2} \over A_{m}} - C_{m} ) }
\prod_{m=1}^{M} \sqrt{ {\pi \over A_{m} } }
\end{eqnarray}

It is instructive to 
make some simplifications.
Assume the variances
are constants, independent of $m$.
(If this is not true, not much simplification
is possible, but the general expressions 
are readily evaluated numerically.)
Then 
\begin{equation}
\sum_{m=1}^{M}
C_{m} =
{1 \over 2}[ { \sum_{m=1}^{M} x_{m}^{2} \over (\sigma^{X})^{2} }
 + { \sum_{m=1}^{M} y_{m}^{2} \over  (\sigma^{Y})^{2} } ]
 \equiv C_{0}
\end{equation}
is just a constant, independent of $\tau$ (and $a$).
Furthermore, 
\begin{equation}
A_{m} \equiv
{1 \over 2} [ { 1 \over (\sigma^{X})^{2} }
 \ + { a^{2} \over (\sigma^{Y})^{2} } ] \equiv A(a)
\end{equation}
is also constant as far as $m$ and $\tau$ are concerned,
although it does depend on $a$.

Hence
\begin{equation}
P( \tau, a | D ) = P_{0}
[{\pi \over A(a) }]^{{M \over 2}}
e^{ - M C_{0} }
e^{ B(\tau, a) }
\end{equation}
where
\begin{eqnarray}
B(\tau, a) & \equiv & 
{\sum_{m=1}^{M} B_{m}^{2} \over A(a)} \\
 & = & {1 \over A(a)} \sum_{m=1}^{M} {1 \over 4}
 [{ x_{m}  \over  (\sigma_{m}^{X})^{2} } +
{ a y_{m+\tau} \over  (\sigma_{m+\tau}^{Y})^{2} } ]^{2} \\
& = & 
 { \sum_{m=1}^{M} x_{m}^{2} \over 4 A(a) (\sigma^{X})^{4} }
 + {\sum_{m=1}^{M} 2 a x_{m} y_{m+\tau} \over 4 A(a) (\sigma^{X})^{2} (\sigma^{Y})^{2} }
+ { a^{2} \sum_{m=1}^{M} y_{m}^{2} \over 4 A(a) (\sigma^{Y})^{4} } \\
& = & 
K_{0}(a) + K_{1}(a) \gamma_{X,Y} (\tau )
\end{eqnarray}
where
\begin{equation}
K_{0}(a)  \equiv  { \sum_{m=1}^{M} x_{m}^{2} \over 4 A(a) (\sigma^{X})^{4} }
+ { a^{2} \sum_{m=1}^{M} y_{m}^{2} \over 4A(a) (\sigma^{Y})^{4} }
\end{equation}
\noindent
and
\begin{equation}
K_{1}(a)  \equiv  {a \over 2 A(a) (\sigma^{X})^{2} (\sigma^{Y})^{2} }
= {a \over (a \sigma^{X})^{2} + (\sigma^{Y})^{2} }
\end{equation}
\noindent
are independent of $m$ and $\tau$,
and
\begin{equation}
\gamma_{X,Y} (\tau ) = \sum_{m=1}^{M} x_{m} y_{m+\tau}
\end{equation}
\noindent
is the ordinary {\it crosscorrelation function}.
Note that the cross term is the only one where
the $\tau$ dependence disappears due to the summation.
Thus we can write
\begin{eqnarray}
P( \tau, a | D ) & = & P_{0}
[{\pi \over A(a) }]^{{M \over 2}}
e^{ - M C_{0} }
e^{ K_{0}(a)   }
e^{  K_{1}(a) \gamma_{X,Y} (\tau )  }
\label{useful}
\end{eqnarray}
again in the same form as in Eq. (\ref{posterior}).
If $a$ is fixed, we have
\begin{equation}
P( \tau, a | D ) \sim e^{K_{1}(a) \gamma_{X,Y} (\tau )} 
%\sim e^{  {a \gamma_{X,Y} (\tau ) \over (a \sigma^{X})^{2} + (\sigma^{Y})^{2}} }
\end{equation}
or
\begin{equation}
logP( \tau, a | D ) \sim K_{1}(a) \gamma_{X,Y} (\tau ) \sim 
{a \gamma_{X,Y} (\tau ) \over (a \sigma^{X})^{2} + (\sigma^{Y})^{2}} 
\end{equation}
\noindent
Eq. (\ref{useful}) can be used to compute
various quantities related to the lag, the scale factor,
and their variances.

\section{Structure in Time Series}

We turn briefly to a different problem,
namely estimating the signal itself.  This section
is an extension of the {\it Bayesian Blocks} 
(\cite{scargle,scargle_1}) method to the case 
where the measurements have a normal error distribution 
and refer to an extended
% are obtained as a weighted integral over a 
range of the independent variable. 

\begin{figure}[htb]
\hspace{3in}
\vspace{1in}
\par
% \hskip -1.2in
\epsfbox{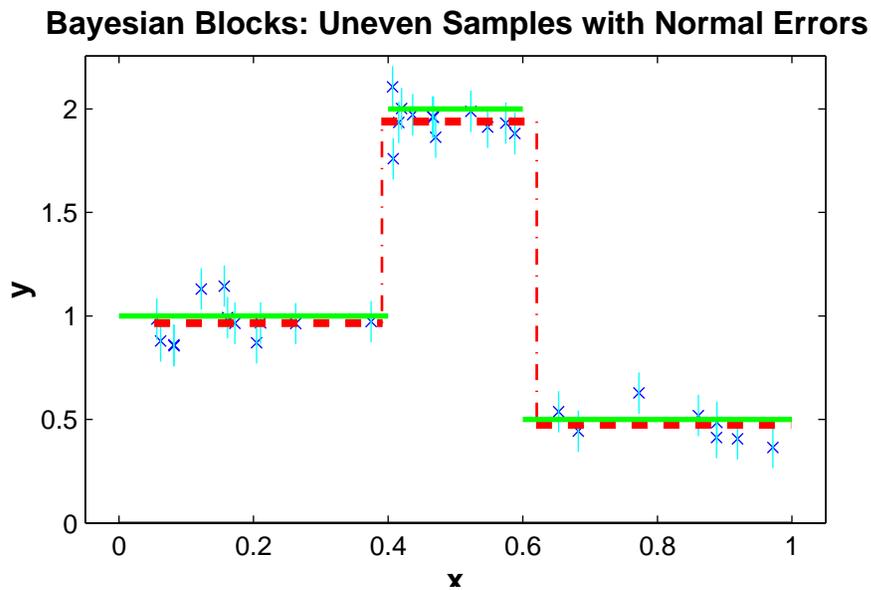}
% \epsfbox{copy1.epsc2}
\caption{Piecewise (block) representations.
Dashed lines: true model. 
Points and error bars: the synthetic data. 
Solid line: Bayesian Block estimate. 
}
\end{figure}

Figure 1 shows the block representation
for a toy problem with just three blocks,
and delta function spread functions for the
independent variable. Note that the 
change point locations have been determined
essentially exactly.

\subsection{The Data}

The data consists of 
measurements of a function $y(x)$,
not actually confined to the single value of $x$
but instead a weighted averaged over a range of $x$-values
\large\begin{equation}
\mbox{Data} = \{ y_{n},  x_{n}, \sigma_{n}, w_{n}(x), n = 1, 2, \dots, N\}
\end{equation}
\noindent
where $\sigma_{n}$ is the known variance of measurement $n$,
and $w_{n}(x)$ is the weighting function,
allowed to be different for each datum.

\subsection{The Model}

We assume the standard piece-wise constant model of the underlying signal,
that is, a set of contiguous blocks:
\begin{equation}
B(x) = \sum_{j=1}^{N_{b}} B^{(j)}( x )
\end{equation}
where each block is represented as a {\it boxcar} function:
\large\begin{eqnarray}
 B^{(k)}( x ) = \{
             \begin{array}{ll}
             B_{j} & \zeta_{j} \le x  \le  \zeta_{j+1}\cr
             0     & \mbox{otherwise}
             \end{array}
\end{eqnarray}
\noindent
the $\zeta_{j}$ are the changepoints, satisfying
\large\begin{equation}
min( x_{n} ) \le \zeta_{1} \le \zeta_{2} \le \dots
   \zeta_{j} \le \zeta_{j+1} \le \dots \le \zeta_{N_{b}} \le max( x_{n} ) 
\end{equation}
\noindent
and the $B_{j}$ are the heights of the blocks.

The value of the observed quantity, $y_{n}$, at $x_{n}$, 
under this model is
\large\begin{eqnarray}
\begin{array}{ll}
\hat{y}_{n}  & = \int w_{n}(x) B( x ) dx \cr
 \           & = \int w_{n}(x) \sum_{j=1}^{N_{b}} B^{(j)}( x ) dx \cr
 \           & = \sum_{j=1}^{N_{b}} \int w_{n}(x) B^{(j)}( x ) dx \cr
 \           & = \sum_{j=1}^{N_{b}} B_{j} \int_{\zeta_{j}}^{\zeta_{j+1}} w_{n}(x) dx 
\end{array}
\end{eqnarray}
\noindent
so we can write
\large\begin{equation}
\hat{y}_{n} = \sum_{j=1}^{N_{b}} B_{j} G_{j}(n) 
\end{equation}
\noindent
where
\large\begin{equation}
G_{j}(n) \equiv \int_{\zeta_{j}}^{\zeta_{j+1}} w_{n}(x) dx 
\end{equation}
\noindent
is the inner product of the $n$-th weight function with the 
support of the $j$-th block.  The analysis in (\cite{bretthorst_1}) 
showns how do deal with the non-orthogonality that 
is generally the case here.\footnote{If 
the weighting functions are delta functions,
it is easy to see that $G_{j}(n)$ is non-zero
if and only if $x_{n}$ lies in block $j$,
and since the blocks do not overlap the product 
$G_{j}(n) G_{k}(n)$ is zero for $j \ne k$,
yielding  orthogonality,
$\sum _{N} G_{j}(n) G_{k}(n) = \delta_{j,k}$.
And of course there can be 
some orthogonal blocks, for which there 
happens to be no``spill over'', but these
are exceptions.}

\subsection{The Posterior}

The averaging process in this data model 
induces dependence among the blocks.
The likelihood, written as a product of
likelihoods of the assumed independent data samples,
is
\begin{eqnarray}
P( \mbox{Data} | \mbox{Model} ) & = \prod_{n=1}^{N} P( y_{n} | \mbox{Model} ) \\
                                & = \prod_{n=1}^{N} {1 \over \sqrt{ 2 \pi \sigma_{n} }}
e^{- {1 \over 2} ({ y_{n} - \hat{y}_{n} \over \sigma_{n} } ) ^{2} } \\
                                & = \prod_{n=1}^{N} {1 \over \sqrt{ 2 \pi \sigma_{n} }} 
e^{- {1 \over 2} ({ y_{n} - \sum_{j=1}^{N_{b}} B_{j} G_{j}(n)  \over \sigma_{n} } ) ^{2} } \\
           & = Q 
e^{- {1 \over 2} ({ y_{n} - \sum_{j=1}^{N_{b}} B_{j} G_{j}(n)  
\over \sigma_{n} } ) ^{2} } \ ,
\end{eqnarray}
\noindent
where
\begin{equation}
Q \equiv \prod_{n=1}^{N} {1 \over \sqrt{ 2 \pi \sigma_{n} }} \ .
\end{equation}
After more algebra
and adopting a new notation, symbolized by
\begin{equation}
{ y_{n} \over \sigma_{n} } \rightarrow y_{n} 
\end{equation}
and
\begin{equation}
{ G_{k}(n) \over \sigma_{n} } \rightarrow G_{k}(n) \ , 
\end{equation}
\noindent
we arrive at
\begin{equation}
log P( \{y_{n}\} | B ) = Q e^{- {H \over 2}} \ ,
\end{equation}
\noindent
where
\large\begin{equation}
H \equiv
\sum_{n=1}^{N} y_{n}^{2}
- 2 \sum_{j=1}^{N_{b}}  B_{j} \sum_{n=1}^{N} y_{n} G_{j}(n) 
+ \sum_{j=1}^{N_{b}}  \sum_{k=1}^{N_{b}}
 B_{j}  B_{k}
\sum_{n=1}^{N}
  G_{j}(n)  G_{k}(n) \ .
\end{equation}
\noindent
The last two equations
are equivalent to Eqs. (3.2) and (3.3) of (\cite{bretthorst_1}),
so that the orthogonalization of the basis functions
and the final expressions follow exactly as in that reference.

%=================================================

\section{Conclusions and Future Work}

This paper has developed an algorithm
for estimating time series lags, leading to
a posterior that is the exponential of a scaled
cross correlation function.
In addition we developed an extension of the
Bayesian Blocks algorithm to the
case where not only are there errors in the
dependent variable, but where the measurement
is a weighted integral over a finite range of the independent 
variable.  
Work planned includes development of numerical algorithms,
testing them on synthetic and real data, 
and then making them freely available 
in the form of Matlab programs.
I also am working on a similar analysis of
scaling behavior in time series, where the
{\it scalegram} -- the square of the wavelet
coefficients, averaged over their location index -- is the sufficient
statistic.

\vskip 0.5 in
% \section{Acknowledgements}
\noindent
{\bf Acknowledgements}

I am greatly indebted to many colleagues
for comments, suggestions, and encouragement:
especially Tom Loredo, Alanna Connors, Larry Bretthorst,
and Peter Sturrock,  
as well as Jay Norris and Jerry Bonnell in connection
with our joint work on Gamma Ray Bursts and the Gamma
Ray Large Area Space Telescope (GLAST).

%=========================================
{99}

\end{document}